\documentclass[english,a4paper,10pt]{article}

\usepackage{xcolor}
\usepackage{amsfonts,amsmath,amstext,amsbsy,amsthm,amscd}
\usepackage{graphics}
\usepackage{epsfig}
\usepackage{geometry}
\usepackage{amssymb}
\usepackage[latin1]{inputenc}
\usepackage{babel,indentfirst,amsfonts,array}

\newtheorem{thm}{Theorem}[section]
\newtheorem{lemme}[thm]{Lemma}

\newtheorem{defi}[thm]{Definition}

\def\T{\mathbb T}
\def\R{\mathbb R}

\def\Q{\mathbb Q}
\def\Z{\mathbb Z}
\def\N{\mathbb N}
\def\({\left(}
\def\){\right)}
\def\[{\left[}
\def\]{\right]}
\def\fin{\hfill\square}

\def\fin{\hfill $\square$}

\title{\textsc{Random self-similar series over a rotation}
\author{Julien Br\'emont}
\date{Universit\'e Paris-Est Cr\'eteil,~janvier 2024}
}

\begin{document}

\maketitle

\setcounter{page}{1}

\begin{abstract}
We study the law of random self-similar series defined above an irrational rotation on the Circle. This provides a natural class of continuous singular non-Rajchman measures.
\end{abstract}

\footnote{
\begin{tabular}{l}\textit{AMS $2020$ subject classifications~: 37E10, 42A38.} \\
\textit{Key words and phrases~: self-similar measure, Rajchman measure, irrational rotation on the torus} 
\end{tabular}}

\section{Introduction}

{\it Dynamical setting.} Consider a probability space $(\Omega,{\cal F},P)$, with a measurable transformation $T:\Omega\rightarrow\Omega$, preserving $P$. The dynamical system $(\Omega,{\cal F},P,T)$ is supposed to be ergodic. 

\medskip
\noindent
Given real random variables $b(\omega)$ and $r(\omega)>0$ on $(\Omega,\cal F)$, define for $\omega\in \Omega$ the real affine map $\varphi_{\omega}(y)=b(\omega)+r(\omega)y$, $y\in\R$. We assume that $\{\varphi_{\omega},~\omega\in\Omega\}=S$ is countable (with $\forall\varphi\in S$, $P(\varphi_{\omega}=\varphi)>0$), $b\in L^1$, $\log r\in L^1$ and $\int_{\Omega} \log r~dP<0$. Setting $r_n(\omega)=r(\omega)\cdots r(T^{n-1}\omega)$, with $r_0(\omega)=1$, introduce the a.-e. defined random variable~:

$$X(\omega)=\sum_{n\geq0}r_n(\omega)b(T^n\omega).$$

\noindent
The law, or occupation measure, of $X$ on $(\R,{\cal B}(\R))$ is denoted by $P_X$, i.e. $P_X(A)=P(X^{-1}(A))$, $A\in {\cal B}(\R)$. The ``self-similar" relation $X(\omega)=\varphi_{\omega}(X(T\omega))$, equivalently rewritten in the ``coboundary" form $b(\omega)=X(\omega)-r(\omega)X(T\omega)$, will be central. It differs from the usual relations of self-similarity for measures, which require some form of independence, not supposed here. Note that if $b(\omega)=\alpha(\omega)-r(\omega)\alpha(T\omega)$, for some random $\alpha$, then necessarily $\alpha=X$, a.-e.. 

\medskip
\noindent
Such a setting includes the traditional self-similar measures (cf Varj\'u \cite{varju} for a survey), corresponding to the independent case, i.e. $\Omega$ a product space with the left shift $T$, $P$ a product measure and $b$, $r$ functions of the first coordinate. Bernoulli convolutions are a famous example, cf the review of Solomyak \cite{solomyak}. The present ergodic extension can be motivated by the case when all affine maps are strict contractions. There is then a self-similar set associated with $S$ and this broader class of measures, supported by $S$, may help studying its properties.

\medskip
A fundamental question concerns the type of $P_X$ with respect to Lebesgue measure $Leb$ and, first of all, the purity of the Radon-Nikodym decomposition. The law of pure types of Jessen and Wintner may be applied to some extent (cf Jessen and Wintner {\cite{JW}, Theorem 35, or Elliott \cite{elliott}, Lemma 1.22), but it seems clearer to give a direct proof in the present situation.

\begin{lemme} The law $P_X$ is of pure type.
\end{lemme}

\medskip
\noindent
{\it Proof of the lemma~:}

\noindent
Let $S^{(n)}=S\circ\cdots\circ S$, $n\geq0$, and $C=\{\varphi\in\cup_{n\geq1}S^{(n)}\mbox{, strict contraction}\}$, countable. Each $\varphi\in C$ having a unique fixed point $fix(\varphi)$, the set ${\cal P}=\{fix(\varphi),~\varphi\in C\}$ is countable.

\medskip
\noindent
- If there exists $a\in\R$, $A=\{X=a\}$, with $P(A)>0$, then $\omega$ a.-e. on $A$, there exists $n\geq 1$ such that $T^n\omega\in A$ and $\varphi_{\omega}\cdots\varphi_{T^{n-1}\omega}\in C$. As $X(\omega)=X(T^n\omega)=a$, we get $a=\varphi_{\omega}\cdots\varphi_{T^{n-1}\omega}(a)$, so $a\in{\cal P}$. Now, $\omega$ a.-e. on $\Omega$, there exists $n\geq0$ such that $T^n\omega \in A$, thus $X(\omega)\in \{\varphi(c),~c\in{\cal P}, \varphi\in \cup_{n\geq0}S^{(n)}\}=:{\cal Q}$, a countable set. Therefore $P_X({\cal Q})=1$ and $P_X$ is purely atomic.

\medskip
\noindent
- If $P_X$ is continuous and if there exists $A\in{\cal B}(\R)$ with $Leb(A)=0$ and $P_X(A)>0$, introduce $B=\cup_{\varphi\in \cup_{n\geq0}S^{(n)}}\varphi^{-1}(A)$. Clearly $Leb(B)=0$. Since $X(\omega)\in B$ implies $X(T\omega)=\varphi_{\omega}^{-1}(X(\omega))\in B$, the set $X^{-1}(B)$ is $T$-invariant. As $P(X^{-1}(B))\geq P(X^{-1}(A))>0$, ergodicity implies that $P_X(B)=P(X^{-1}(B))=1$. Therefore $P_X\perp Leb$. \fin

\bigskip
{\it Pure atomicity.} Let us discuss the continuity of $P_X$. Clearly, $P_X=\delta_c$ if and only if $\forall \varphi\in S$, $\varphi(c)=c$. In the independent case, the purely atomic situation reduces to $P_X$ a Dirac mass, as follows from the relation (obtained when conditioning with respect to the first step)~:

$$P_X(A)=\sum_{\varphi\in S}P(\varphi_{\omega}=\varphi)P_X(\varphi^{-1}(A)),~A\in{\cal B}(\R).$$

\noindent
Indeed, if there exists an atom, then the latter implies that the non-empty finite set $E$ of points defining an atom of maximal mass is stable under any $\varphi^{-1}$. Finiteness of an orbit under iterations of an affine map forces any $c\in E$ to be a fixed point of any $\varphi\in S$. 

\medskip
\noindent
This is far from true in the general ergodic context. Fixing $r$ and any $\alpha\in L^1$ with countable support, when setting $b=\alpha-r\alpha\circ T$, we have $X=\alpha$. As a result, $P_X$ can be discrete with even non-finite support. Moreover, as we shall see later, determining the conditions under which $P_X$ is continuous can be a non-degenerate problem.

\medskip
\noindent
Mention here a recipe for building non-trivial examples of discrete laws when $r(\omega)=\lambda\in(0,1)$ is algebraic. Let for instance $\lambda=0,618...$ be the inverse of the Golden Mean, i.e. $\lambda^2+\lambda-1$=0. Taking $g\in L^1$ with countable support and $b=g+g\circ T-g\circ T^2$, then $b=(g+(1+\lambda)g\circ T)-\lambda(g\circ T+(1+\lambda)g\circ T^2)$. This means that $X(\omega)=g(\omega)+(1+\lambda)g(T\omega)$. 

\medskip
\noindent
More generally, if $\sum_{k=0}^p\alpha_k\lambda^{p-k}=0$, $p\geq1$, let $b(\omega)=\sum_{k=0}^p\alpha_kg(T^k\omega)$, where $g\in L^1$ has countable support. Then $X(\omega)=\sum_{n=0}^{p-1}g(T^n\omega)(\sum_{k=0}^n\alpha_k\lambda^{n-k})$, as  $X(\omega)-r(\omega)X(T\omega)=b(\omega)$.

\medskip
Recall also the link between the existence of atoms and the Fourier transform. We define~:

$$\hat{P}_X(t)=\int_{\R}e^{2i\pi tx}dP_X(x),~t\in\R.$$

\smallskip
\noindent
If $P_X$ is continuous, then, by Wiener's theorem~:

$$\frac{1}{R}\int_0^R|\hat{\mu}(t)|^2dt\rightarrow0\mbox{, as }R\rightarrow+\infty.$$

\smallskip
\noindent
A more precise information of local regularity is when $P_X$ is a Rajchman measure, meaning that $\hat{P}_X(t)\rightarrow0$, as $t\rightarrow+\infty$. Equivalently, $tX\mod 1\rightarrow_{\cal L}Leb_{\T}$, as $t\rightarrow+\infty$. A classical example of continuous non-Rajchman measures is the uniform measure on the triadic Cantor set. The present paper furnishes a natural class of such measures.




\bigskip
{\it Content of the article.} We study the special case when the dynamics is given by an irrational rotation on the 1-torus, with functions $b$ and $r$ locally constant on some finite collection of intervals. For obvious complexity reasons, $P_X$ is singular, even of zero-dimensional support, so it remains to decide between continuous singularity and pure atomicity. We show that the latter is equivalent to the simultaneous satisfaction of a finite number of explicit algebraic equations. Generically, $P_X$ appears to be continuous, but also not a Rajchman measure. In the last section, we discuss another approach of the continuity problem for general systems.

\section{The case of the Circle}

Let $\T=\R\backslash\Z$ be the 1-torus, with uniform measure $Leb_{\T}$ and an irrational rotation $T$ of angle $\alpha\in(0,1)$. We recall classical material about continued fractions; see for example Khinchin's book \cite{khinchin_book}. The angle $\alpha$ can be expanded in infinite continued fraction~:

$$\alpha=\cfrac{1}{a_1+\cfrac{1}{a_2+\cdots}}=[0,a_1,a_2,\cdots],$$

\medskip
\noindent
where the partial quotients $(a_i)_{i\geq1}$ are obtained by iterations of the Gauss map, starting from $\alpha$. The successive truncations $[0,a_1,a_2,\cdots,a_n]=p_n/q_n$, $n\geq1$, are the convergents of $\alpha$. The $(p_n)$ and $(q_n)$ check the same recursive relation:

$$p_{n+1}=a_{n+1}p_n+p_{n-1},~~~q_{n+1}=a_{n+1}q_n+q_{n-1},~~n\geq0,$$  

\medskip
\noindent
with $p_0=0,p_{-1}=1$ and $q_0=1,q_{-1}=0$. Classical inequalities are (cf \cite{khinchin_book}, chap. 1)~:

$$\frac{1}{2q_{n+1}}\leq\frac{1}{q_n+q_{n+1}}\leq \Vert q_n\alpha\Vert\leq\frac{1}{q_{n+1}},$$

\medskip
\noindent
where $\|x\|$ is the distance from $x$ to $\Z$. Our purpose is to establish the following result.

\begin{thm} 

$ $

\noindent
Let $T$ be a rotation of angle $\alpha=[0,a_1,a_2,\cdots]\not\in\Q$ on $\T$. 

\medskip
\noindent
Given $N\geq1$ points $d_0<d_1<\cdots<d_{N-1}<d_N=d_0$ on $\T$, consider on ${\cal D}=\{d_0,\cdots,d_{N-1}\}$ the partial order ``$d_i\rightarrow d_j$ iff $d_j=T^pd_i$ for some $p\geq 0$". Partition ${\cal D}=\sqcup_{1\leq k\leq K}{\cal D}_k$ into maximal subsets ${\cal D}_k=\{d_{0,k}\rightarrow\cdots\rightarrow d_{m_k,k}\}$, with $m_k\geq0$; define $p_k\geq0$ by $d_{m_k,k}=T^{p_k}d_{0,k}$.

\medskip
\noindent
Let $b:\T\rightarrow\R$ and $r=\T\rightarrow(0,1)$ be constant on each interval $[d_i,d_{i+1})$, $0\leq i<N$. Define $X(x)=\sum_{n\geq 0}b(T^nx)r_n(x)$, $x\in\T$, and denote by $P_X$ the image of $Leb_{\T}$ by $X$. Then~:

\begin{enumerate}
\item $Supp(P_X)$ has box-counting dimension zero, in particular $P_X\perp Leb$.
\item The measure $P_X$ is continuous iff $X$ is discontinuous at some $d_{0,k}$, $1\leq k\leq K$. Otherwise $X$ is constant on the intervals of the partition determined by $\{T^pd_{0,k},~0\leq p\leq p_k,~0\leq k\leq K\}$, hence $\mbox{Supp}(P_X)$ is finite, with at most $\sum_{1\leq k\leq K}(1+p_k)$ elements.
\item If $a_n\geq 10N+20N^2\ln13/(-\ln\|r\|_{\infty})$ infinitely often, then $P_X$ is not a Rajchman measure. If $(a_n)$ is unbounded, then $t_nX\mod 1\rightarrow_{{\cal L}}0$, along a sequence of integers $(t_n)\rightarrow+\infty$ . 
\end{enumerate}
\end{thm}

\medskip
\noindent
{\it Proof of the theorem~:}

\noindent
$1)$ For any $n\geq1$, $x\longmapsto\sum_{k=0}^{n-1}r_k(x)b(T^kx)$ is constant on each interval of the partition determined by $\cup_{0\leq k<n}T^{-k}{\cal D}$, so takes at most $nN$ values. As $|\sum_{k\geq n}r_k(x)b(T^kx)|\leq \|r\|_{\infty}^n\|b\|_{\infty}/(1-\|r\|_{\infty})$, $Supp(P_X)$ can be covered for any $\varepsilon>0$ by at most $-C\log \varepsilon$ balls of radius $\varepsilon$, for some constant $C>0$. This gives the result.

\medskip
\noindent
$2)$ In the present context of strict contractions, $X$ is right-continuous and admits a left limit $X(x^-)$ at every $x\in\T$. Set $\Delta_k=X(d_{0,k})-X(d_{0,k}^-)$ and ${\cal K}=\{1\leq k\leq K,~\Delta_k\not=0\}$. Supposing that ${\cal K}\not=\o$, we set $\Delta=\min_{k\in{\cal K}}|\Delta_k|>0$. Choose also $\varepsilon>0$ so that~:

\begin{equation}
\label{epsi}
\min_{k\in {\cal K}}\inf_{\underset{|y-x|\leq \varepsilon}{x<d_{0,k}\leq y}}|X(x)-X(y)|\geq\Delta/2.
\end{equation}

\noindent
Set $\rho_k^{\pm}=r_{p_k+1}(d_{0,k}^{\pm})$, $1\leq k\leq K$, and define $\rho_{\max/\min}=\max/\min\{\rho_k^{\pm},~1\leq k\leq K\}$. For the sequel, fix $M>\max\{p_1,\cdots,p_K\}$ such that~:

\begin{equation}
\label{M}
\|X\|_{\infty}\|r\|_{\infty}^{M-1}<\frac{\Delta}{12N}\({\frac{\rho_{\min}}{\rho_{\max}}}\)^{3N}.
\end{equation}

\medskip
\noindent
For $1\leq k\leq K$, call $(T^pd_{0,k})_{0\leq p\leq p_k}$ the chain $C_k$. Choose $\gamma(M)>0$ such that for any $x<y<x+\gamma(M)$, each interval $T^k(x,y]$, $k\geq0$, meets at most one element of ${\cal D}$ and after covering the last element of a chain the (necessarily) first element of the next chain is not met until $M$ steps. 

\medskip
\noindent
Take $x\not\in\cup_{l\geq0}T^{-l}{\cal D}$ and $0<\gamma_x<\min\{\gamma(M),\varepsilon\}$ such that if $x<y<x+\gamma_x$, then $T^k(x,y]$ meets no $d_j$, for $0\leq k\leq M$. If $T^k(x,y]$ meets for the first time a chain, it thus has to be at the first element of the chain. For the moment, fix $y$ like this. The choice of $x,y$ is precised later.

\medskip
\noindent 
We consider $X(x)-X(y)$. This way, let $0=t_0<s_1<t_1<s_2<t_2<\cdots$, where, for $i\geq 0$, the $[t_i,s_{i+1})$ are the maximal time intervals of $k$ where $T^k(x,y]$ meets no chain. For $i\geq1$, the $(T^k(x,y])_{k\in [s_i,t_i)}$ cover some chain, say $C_{l_i}$, with $d_{0,l_i}\in T^{s_i}(x,y]$ and $d_{m_{l_i},l_i}\in T^{t_i-1}(x,y]$. 

\medskip
\noindent
Introduce $r_n(x)=r_{s_n-t_{n-1}}(T^{t_{n-1}}x)$, $n\geq1$. We define $n_0\geq1$ as the first integer $n$ such that $l_n\in{\cal K}$. First of all, we can write~:

$$X(x)-X(y)=r_1(x)(X(T^{s_1}x)-X(T^{s_1}y)).$$

\medskip
\noindent
In a recursion, suppose now that for some $1\leq n<n_0$~:

\begin{equation}
\label{formula}
X(x)-X(y)=r_1(x)\cdots r_{n}(x)\sum_{0\leq u<n}\rho_1^*\cdots \rho_{n-1}^*(X(x_{u}^n)-X(x_{u+1}^n)),
\end{equation}

\medskip
\noindent
with points $T^{s_n}x=x_0^n\leq x_1^n\leq \cdots\leq x_{n}^n=T^{s_n}y$ and $\rho_i^*=\rho_{l_i}^{\pm}$. Since $T^{s_n}x<d_{0,l_n}\leq T^{s_n}y$, let $v$ be the index such that $x_{v}^n<d_{0,l_n}\leq x_{v+1}^n$. Adding $d_{0,l_n}$ to the $(x_{i}^n)_{0\leq i\leq n}$ gives $n+2$ points, written in their natural order as $(y_u^n)_{0\leq u\leq n+1}$. Since $n<n_0$, we split in the following way the term for $u=v$ in \eqref{formula}~:

\begin{eqnarray}
X(x_{v}^n)-X(x_{v+1}^n)
&=&X(x_{v}^n)-X(d_{0,l_n}^-)+X(d_{0,l_n})-X(x_{v+1}^n)\nonumber\\
&=&X(y_{v}^n)-X(y_{v+1}^{n,-})+X(y_{v+1}^n)-X(y_{v+2}^n).\nonumber
\end{eqnarray}

\medskip
\noindent
Set $\rho_n^*=\rho_{l_n}^-$ if $u\leq v$ and $\rho_n^*=\rho_{l_n}^+$ if $u\geq v+1$. For $u\not=v$~:

\begin{eqnarray}
X(y_{u}^n)-X(y_{u+1}^n)&=&\rho_{l_n}^*(X(T^{t_n-s_n}y_{u}^n)-X(T^{t_n-s_n}y_{u+1}^n))\nonumber\\
&=&\rho_{l_n}^*r_{n+1}(x)(X(T^{s_{n+1}-s_n}y_{u}^n)-X(T^{s_{n+1}-s_n}y_{u+1}^n)).\nonumber
\end{eqnarray}

\medskip
\noindent
Now, in the same way~:

$$X(y_{v}^n)-X(y_{v+1}^{n,-})=\rho_{l_n}^-r_{n+1}(x)(X(T^{s_{n+1}-s_n}y_{v}^n)-X(T^{s_{n+1}-s_n}y_{v+1}^{n,-})).$$

\medskip
\noindent
As $T^{s_{n+1}-s_n}y^n_{v+1}=T^{s_{n+1}-t_n+1}d_{m_{l_n},l_n}$ and $s_{n+1}-t_n+1\geq1$, from the continuity of $X$ at any $T^kd_{m_{l_n},l_n}$, $k\geq1$, we get $X(T^{s_{n+1}-s_n}y_{v+1}^{n,-})=X(T^{s_{n+1}-s_n}y_{v+1}^{n})$. We can now finally set $x_{u}^{n+1}=T^{s_{n+1}-s_n}y_{u}^n$, $0\leq u\leq n+1$, and we obtain when replacing in \eqref{formula} that the latter is satisfied with $n$ replaced by $n+1$. As a result, the formula is true for $n=n_0$~:

\begin{equation}
\label{formula2}
X(x)-X(y)=r_1(x)\cdots r_{n_0}(x)\[{\sum_{0\leq u<n_0}\rho_1^*\cdots \rho_{n_0-1}^*(X(x_{u})-X(x_{u+1}))}\],
\end{equation}

\medskip
\noindent
with, simplifying notations, points $T^{s_{n_0}}x=x_0\leq x_1\leq \cdots\leq x_{n_0}=T^{s_{n_0}}y$ and $\rho_i^*=\rho_{l_i}^{\pm}$. Again $T^{s_{n_0}}x<d_{0,l_{n_0}}\leq T^{s_{n_0}}y$ and let $v$ be the index such that $x_{v}<d_{0,l_{n_0}}\leq x_{v+1}$.

\medskip
\noindent
Now, using \eqref{epsi}, by definition, $|X(x_{v})-X(x_{v+1})|\geq\Delta/2$, whereas, as before, for $u\not =v$~:

$$X(x_{u})-X(x_{u+1})=\rho_{l_{n_0}}^*r_{n_0+1}(x)(X(T^{s_{n_0+1}-s_{n_0}}x_{u})-X(T^{s_{n_0+1}-s_{n_0}}x_{u+1})).$$

\medskip
\noindent
Since $M$ verifies $r_{n_0+1}(x)=r_{s_{n_0+1}-t_{n_0}}(T^{t_{n_0}}x)\leq \|r\|_{\infty}^{M-1}$, when calling $A$ the term between brackets in \eqref{formula2}, we deduce from the previous considerations that~:

\begin{eqnarray}
|A|&\geq& \frac{\Delta}{2}(\rho_{\min})^{n_0-1}-2\|X\|_{\infty}(n_0-1)(\rho_{\max})^{n_0}r_{n_0+1}(x)\nonumber\\
\label{A}&\geq&\frac{(\rho_{\min})^{n_0}}{2}\[{\Delta-4n_0\|X\|_{\infty}\({\frac{\rho_{\max}}{\rho_{\min}}}\)^{n_0}\|r\|_{\infty}^{M-1}}\].
\end{eqnarray}

\medskip
Suppose $P_X$ purely atomic. Let $x$ be a Lebesgue density point in some atom ($Leb_{\T}$ a.-e. point is such a point), not in the countable set $\cup_{l\geq0}T^{-l}{\cal D}$. Choose $n$ large enough so that $3\|q_n\alpha\|<\gamma_x$ and take $y\in x+(2\|q_n\alpha\|,3\|q_n\alpha\|)$ verifying $X(x)=X(y)$. This is possible, as the proportion of points in $x+(0,3\|q_n\alpha\|)$ lying in the same atom as $x$ tends to one, as $n\rightarrow+\infty$.

\medskip
\noindent
Recall that the $(0,\|q_n\alpha\|)+k\alpha$, $0\leq k<q_{n+1}$, are disjoint and, as a classical consequence of the identity $q_n\|q_{n+1}\alpha\|+q_{n+1}\|q_n\alpha\|=1$, that the $x+(0,2\|q_n\alpha\|)+k\alpha$, $0\leq k<q_{n+1}$, cover $\T$, each point belonging to at most two intervals. 

\medskip
\noindent
As a result, the Circle $\T$ is covered by the $T^k(x,y]$, $0\leq k<q_{n+1}$, and each point of $\T$ is covered at most 3 times. We deduce that the $T^k(x,y]$ will pass at most three times in chains $C_z$, $z\not\in {\cal K}$, before finally meeting a chain whose index is in ${\cal K}$. Therefore $n_0\leq 3N$. From \eqref{A}~:

$$|A|\geq\frac{(\rho_{\min})^{n_0}}{2}\[{\Delta-12N\|X\|_{\infty}\({\frac{\rho_{\max}}{\rho_{\min}}}\)^{3N}\|r\|_{\infty}^{M-1}}\]>0,$$

\noindent
using property \eqref{M} of $M$. Since $A\not=0$ and $r_1(x)\cdots r_{n_0}(x)\not=0$, we get a contradiction in \eqref{formula2} with the fact that $X(x)-X(y)=0$.

\medskip
In the other direction, suppose that $\Delta_k=0$, $1\leq k\leq K$. The set $\{T^pd_{0,k},~0\leq p\leq p_k,~1\leq k\leq K\}$, the union of the chains, gives a partition of $\T$ into $\sum_{k=1}^K(1+p_k)$ intervals. We show that $X$ is constant on each piece. This way, let $M>2+\max\{p_1,\cdots,p_k\}$ and take the corresponding $\gamma(M)>0$. Take $x<y$ interior to the same interval of the partition, with $x<y<x+\gamma(M)$. Considering the orbit $T^k(x,y]$, $k\geq0$, if a chain is met for the first time, then it is at the first element of the chain. As $\Delta_k=0$ for all $1\leq k\leq K$, formula \eqref{formula} is true for all $n\geq1$~:

$$X(x)-X(y)=r_1(x)\cdots r_{n}(x)\sum_{0\leq u<n}\rho_1^*\cdots \rho_{n-1}^*(X(x_{u}^n)-X(x_{u+1}^n)),$$

\noindent
with, using the same notations for time intervals, points $T^{s_n}x=x_0^n\leq x_1^n\leq \cdots\leq x_{n}^n=T^{s_n}y$ and $\rho_i^*=\rho_{l_i}^{\pm}$. As $r_k(x)\leq \|r\|_{\infty}^{M-1}\leq \|r\|_{\infty}$, we get~:

$$|X(x)-X(y)|\leq\|r\|_{\infty}^n\times n\rho_{\max}^{n-1}\times2\|X\|_{\infty}.$$

\smallskip
\noindent
As this goes to 0, as $n\rightarrow+\infty$, we get $X(x)=X(y)$. Hence $X$ is locally constant, hence constant, on each interval of the partition. This concludes the proof of point $2)$.

\bigskip
\noindent
$3)$ We examine the Rajchman character of $P_X$. Set $S_k(x)=-\sum_{l=0}^{k-1}\log r(T^lx)$, with $S_0=0$. Then $X(x)=\sum_{k\geq0}e^{-S_k(x)}b(T^kx)$. Fixing $n$ and $0\leq m_n\leq a_{n+1}$, arbitrary for the moment~:

\begin{eqnarray}
X(x)&=&\sum_{k=0}^{q_n-1}e^{-S_k(x)}\sum_{m\geq0}e^{-S_{mq_n}(T^kx)}b(T^{mq_n+k}x)\nonumber\\
\label{un}&=&\sum_{k=0}^{q_n-1}e^{-S_k(x)}\sum_{0\leq m\leq m_n}e^{-S_{mq_n}(T^kx)}b(T^{mq_n+k}x)\\
\label{trois}&+&\sum_{k=0}^{q_n-1}e^{-S_k(x)}\sum_{m>m_n}e^{-S_{mq_n}(T^kx)}b(T^{mq_n+k}x).\end{eqnarray}

\medskip
\noindent
Suppose $n$ even (the other case is similar), so $q_n\alpha\mod 1$ is on the right side of 0 on the Circle. Consider  \eqref{un} and $0\leq k<q_n$, as well as $m\geq1$. If $[T^{k+l}x,T^{k+l+(m-1)q_n}x]$ contains no $d_i$, for any $0\leq l<q_n$, then $S_{mq_n}(T^kx)=mS_{q_n}(T^kx)$. Similarly, $b(T^{mq_n+k}x)=b(T^kx)$, whenever $[T^{k}x,T^{k+mq_n}x]$ contains no $d_i$. Introduce~:

$$\Omega_n=\underset{0\leq k<2q_n,0\leq i<N}{\cup}-k\alpha-d_i+[-m_nq_n\alpha,0],$$ 

\medskip
\noindent
of measure $\leq 2q_nNm_n\|q_n\alpha\|\leq 2Nm_n/a_{n+1}$. For $x\not\in\Omega_n$, one has $X(x)=Z_n(x)+R_n(x)$, with~:

\begin{eqnarray}
Z_n(x)=\sum_{k=0}^{q_n-1}e^{-S_k(x)}b(T^kx)\frac{1-e^{-(m_n+1)S_{q_n}(T^kx)}}{1-e^{-S_{q_n}(T^kx)}},~\|R_n\|_{\infty}\leq \frac{\|b\|_{\infty}\|r\|_{\infty}^{(m_n+1)q_n}}{1-\|r\|_{\infty}}.\nonumber
\end{eqnarray}

\medskip
\noindent
For any $t_n>0$, decomposing $e^{2i\pi t_n(Z_n+R_n)}-1=e^{2i\pi t_nZ_n}(e^{2i\pi t_nR_n}-1)+e^{2i\pi t_nZ_n}-1$ and using that $x\longmapsto e^{ix}$ is 1-Lipschitz on $\R$, we have~:

\begin{eqnarray}
\label{deco}
|\hat{P}_X(t_n)-1|&\leq&\int_{\Omega_n^c}|e^{2i\pi t_nX}-1|dx+2|\Omega_n|\\
&\leq&\int_{\Omega_n^c}|e^{2i\pi t_nZ_n(x)}-1|dx+t_n\|R_n\|_{\infty}|\Omega_n^c|+4Nm_n/a_{n+1}.\nonumber
\end{eqnarray}

\medskip
\noindent
Now, $Z_n$ is constant on each interval of the partition determined by $\cup_{0\leq l<2q_n}T^{-l}{\cal D}$ and therefore takes at most $2Nq_n$ values. Fixing an integer $r_n\geq4$, cut the torus $\T^{2Nq_n}$ in cubes of sides of length $1/r_n$. This gives $r_n^{2Nq_n}$ cubes. Considering the integers $\{nk,~0\leq k\leq r_n^{2Nq_n}\}$, by the pigeonhole principle, there exists an integer $nt_n$, with $1\leq t_n\leq r_n^{2Nq_n}$, such that $\|nt_nZ_n(x)\|\leq 1/r_n$, for all $x\in\T$. Replacing $t_n$ by $nt_n$ (arbitrary large)~:

\begin{eqnarray}
|\hat{P}_X(nt_n)-1|&\leq&|\Omega_n^c|2\pi/r_n+nt_n\|R_n\|_{\infty}+4Nm_n/a_{n+1}\nonumber\\
&\leq&2\pi/r_n+nr_n^{2Nq_n}\frac{\|b\|_{\infty}\|r\|_{\infty}^{(m_n+1)q_n}}{1-\|r\|_{\infty}}+4Nm_n/a_{n+1}.\nonumber\end{eqnarray}

\noindent
We shall impose $m_n\geq \ln(r_n^{2N})/(-\ln \|r\|_{\infty})$, giving~:

\begin{equation}
\label{pasrajch}
|\hat{P}_X(nt_n)-1|\leq2\pi/r_n+4Nm_n/a_{n+1}+n\|r\|_{\infty}^{q_n}\frac{\|b\|_{\infty}}{1-\|r\|_{\infty}}.
\end{equation}

\noindent
If $r_n\geq 4\pi$ and $m_n\leq a_{n+1}/(10N)$, then $|\hat{P}_X(nt_n)-1|\leq 1/2+2/5+o(1)=9/10+o(1)$. Fixing $r_n=13>4\pi$, then $P_X$ is not a Rajchman measure whenever for infinitely many $n$, one can find an integer $m_n$ satisfying the inequalities~:

\begin{equation}
\label{ineq}
2N\ln r_n/(-\ln\|r\|_{\infty})\leq m_n\leq a_{n+1}/(10N).
\end{equation}

\medskip
\noindent
Since $r_n=13$, this is thus true $a_{n+1}/(10N)\geq1+2N\ln13/(-\ln\|r\|_{\infty})$, along a subsequence.

\medskip
If the partial quotients are unbounded, take~:

$$r_n=a_{n+1}\mbox{ and }m_n=[\sqrt{a_{n+1}}],$$

\noindent
along a subsequence where $a_{n+1}\rightarrow+\infty$. Then \eqref{ineq} is true for large $n$. By \eqref{pasrajch}, $\hat{P}_X(nt_n)\rightarrow1$ along a subsequence $nt_n\rightarrow+\infty$. Next, for any integer $m\geq1$, $|e^{2i\pi t_nmX}-1|\leq m|e^{2i\pi t_nX}-1|$. Keeping the same sequence $(nt_n)$, relation \eqref{deco} at time $nt_n$ for $mX$ gives~:

$$|\hat{P}_X(mnt_n)-1|\leq m\int_{\Omega_n^c}|e^{2i\pi nt_nX}-1|dx+2|\Omega_n|.$$

\medskip
\noindent
As before, the integral and $|\Omega_n|$ go to zero, as $n\rightarrow+\infty$, along the above mentioned subsequence. This completes the proof of point $3)$.

\fin

\bigskip
\noindent
\begin{remark} Explicitly, $P_X$ is purely atomic if and only if for all $1\leq k\leq K$~:

$$\sum_{i=0}^{p_k}\[{r_i(d_{0,k})b(T^id_{0,k})-r_i(d_{0,k}^-)b(T^id_{0,k}^-)}\]+\[{r_{p_k+1}(d_{0,k})-r_{p_k+1}(d_{0,k}^-)}\]X(Td_{m_k,k})=0.$$

\noindent
Because of $X(Td_{m_k,k})$, this value may involve the whole orbit of $d_{0,k}$. On the contrary, when $r(x)=\lambda\in(0,1)$ and writing any maximal set as ${\cal D}_k=\{d_{0,k}\rightarrow_{p_{0,k}}\cdots\rightarrow_{p_{m_k-1,k}} d_{m_k,k}\}$, with integers $p_{i,k}\geq1$ such that $d_{i+1,k}=T^{p_{i,k}}d_{i,k}$, this simplifies into~:

$$\sum_{i=0}^{m_k}\lambda^{p_{0,k}+\cdots+p_{i-1,k}}\[{b(d_{i,k})-b(d_{i,k}^-)}\]=0,~1\leq k\leq K.$$
\end{remark}

\bigskip
\noindent
\begin{remark}
\noindent
If for example all $d_i$ are in distinct orbits, the condition of pure atomicity reduces to $b(d_i)-b(d_i^-)+[r(d_i)-r(d_i^-)]X(Td_i)=0$ and, when $r(x)$ is constant, to $b(d_i)-b(d_i^-)=0$, $0\leq i<N$, i.e. $b$ constant, thus giving $P_X=\delta_{b/(1-\lambda)}$. Proceeding as indicated in the Introduction, it is easy to build examples with any finitely supported law.
\end{remark}

\bigskip
\noindent
\begin{remark}
Concerning point $3)$, we conjecture that $P_X$ is never a Rajchman measure. Here is a classical situation where the result is true for any angle. Recall that a Pisot number $\rho>1$ is an algebraic integer, with Galois conjugates of modulus $<1$. 
\end{remark}

\begin{lemme}

$ $

\noindent
Let $T$ be a rotation of angle $\alpha$ on $\T$, $r(x)=\lambda\in(0,1)$, with $1/\lambda$ a Pisot number, and $b(x)\in\Z$, locally constant on a partition $\T=\sqcup_{0\leq i<N}[d_i,d_{i+1})$. Then $P_X$ is not a Rajchman measure.
\end{lemme}

\noindent
{\it Proof of the lemma~:}

\noindent
In this case, $X(x)=\sum_{k\geq0}\lambda^kb(T^kx)$. If $B\subset\Z$ denotes the finite set of values of $b$, then~:

$$\mbox{Supp}(P_X)\subset\left\{{\sum_{k\geq0}\lambda^kb_k,~b_k\in B}\right\}.$$

\medskip
\noindent
Classically, the latter self-similar set is a set of uniqueness for trigonometric series, hence cannot support a Rajchman measure; cf for example the general result of Varj\'u-Yu \cite{VH}, Theorem 1.4. 

\medskip
For a more elementary proof, introduce the conjugates $\mu_1,\cdots,\mu_d$ of $1/\lambda$ and recall that $\lambda^{-n}+\mu_1^n+\cdots+\mu_d^n\in\Z$, $n\geq0$. If $P_X$ were a Rajchman measure, we would have in particular $\lambda^{-n}X\mod1\rightarrow_{\cal L}Leb_{\T}$, hence $\lambda^{-n}X\circ T^{-n}\mod1\rightarrow_{\cal L}Leb_{\T}$. However, modulo $1$~:

$$\lambda^{-n}X(T^{-n}x)\equiv\sum_{k=1}^n\lambda^{-k}b(T^{-k}x)+X(x)\equiv X(x)-\sum_{k=1}^n(\mu_1^k+\cdots+\mu_d^k)b(T^{-k}x).$$

\medskip
\noindent
The term on the right-hand side converges pointwise to the real random variable~:

$$Y(x)=X(x)-\sum_{k\geq1}(\mu_1^k+\cdots+\mu_d^k)b(T^{-k}x),$$

\noindent
We would get $P_{Y\mod1}=Leb_{\T}$, on $\T$. However, $Y_n(x)\rightarrow Y(x)$, as $n\rightarrow+\infty$, where~:

$$Y_n(x)=\sum_{k=0}^n\lambda^kb(T^kx)-\sum_{k=1}^n(\mu_1^k+\cdots+\mu_d^k)b(T^{-k}x).$$
 
\noindent
We have $\|Y-Y_n\|_{\infty}\leq C\rho^n$, where $\rho=\max\{\lambda,|\mu_1|,\cdots,|\mu_d|\}<1$. Since $Y_n$ takes at most $(2n+1)N$ values, we get $Leb(\mbox{Supp}(P_Y)=0$. Hence $P_Y$ on $\R$ is singular. Therefore $P_{Y\mod1}$ is singular on $\T$ and in particular $P_{Y\mod1}\not=Leb_{\T}$. This concludes the proof of the lemma.

\fin

\section{A remark for general dynamical systems}

For the general setting of the Introduction, we discuss in this last section another approach, relating the continuity of the measure $P_X$ to a question of fixed points. We suppose the dynamical system ergodic and invertible.

\medskip
\noindent
Changing notations, write $\varphi_{\omega}=\psi_{\epsilon(\omega)}$, $\epsilon(\omega)\in {\cal S}$, where ${\cal S}$ is a countable set. For simplicity, we suppose that all affine maps $\psi_j$, $j\in{\cal S}$, are strict contractions. We shall use multi-indices $i=(i_0,\cdots,i_{n-1})\in {\cal S}^n$, for $n\geq1$. We also write $\psi_i=\psi_{i_0}\cdots\psi_{i_{n-1}}$.

\begin{defi}
A multi-index $i\in {\cal S}^n$, $n\geq1$, is minimal if $P((\epsilon,\cdots,T^{n-1}\epsilon)=i)>0$ and for any strict prefix $j$ of $i$, $fix(\psi_j)\not=fix(\psi_i)$. Let ${\cal M}=\{i\in\cup_{n\geq1}{\cal S}^n\mbox{, minimal}\}$.
\end{defi} 
 
\noindent
\begin{remark}
It is easily verified that $fix(\psi_i)=fix(\psi_j)$ if and only if $\psi_i\circ\psi_j=\psi_j\circ\psi_i$. 
\end{remark}

 \begin{lemme}
 
 $ $
 
 \noindent
Suppose the map : $i$ minimal $\longmapsto fix(\psi_i)$, from ${\cal M}$ to $\R$, injective. Then, either $P_X$ is continuous or there exists $N\geq1$ and $(i_0,\cdots,i_{N-1})$ such that for a.-e. $\omega$, $(\epsilon(T^n\omega))_{n\geq0}$ is a left shift of the periodic sequence $(\overline{i_0,\cdots,i_{N-1}},\cdots)\in{\cal S}^{\N}$, in which case $X(\Omega)=\{\psi_{i_k}\cdots\psi_{i_{N-1}}(c),~0\leq k<N\}$, up to a null set, where $c=fix(\psi_{i_0}\cdots\psi_{i_{N-1}})$.
 \end{lemme}
 
\noindent
{\it Proof of the lemma~:}

\noindent
If $P_X$ is purely atomic, let $c$ and $A=\{X=c\}$, with $P_X(A)>0$. On $A$, let $\tau\geq1$ be the return time, a.-e. defined. Then, restricting to sequences appearing with positive probability, $(\epsilon(\omega),\cdots,\epsilon(T^{\tau(\omega)-1}\omega))$ is minimal, as $c=\psi_{\epsilon(\omega)}\cdots\psi_{\epsilon(T^{\tau(\omega)-1}\omega)}(c)$ and if $c=\psi_{\epsilon(\omega)}\cdots\psi_{\epsilon(T^{m-1}\omega)}(c)$ for some $m<\tau(\omega)$, then $X(T^m\omega)=c$, by injectivity, contradicting the definition of $\tau(\omega)$. 

\medskip
\noindent
Since for a.-e. $\omega\in A$, $(\epsilon(\omega),\cdots,\epsilon(T^{\tau(\omega)-1}\omega))$ is minimal and $c$ is the corresponding fixed point, the hypothesis implies that there exists $N\geq1$ and $(i_0,\cdots,i_{N-1})\in{\cal S}^N$ such that $\tau(\omega)=N$ and $(\epsilon(\omega),\cdots,\epsilon(T^{N-1}\omega))=(i_0,\cdots,i_{N-1})$, for a.-e. $\omega$ in $A$. Also, clearly, $X=c$, a.-e. on $A$.

\medskip
\noindent
By ergodicity and invertibility, we now have, up to a null set, $\Omega=\sqcup_{0\leq k<N}T^kA$. Then, for a.-e. $\omega$, the sequence $(\epsilon(T^n\omega))_{n\geq0}$ is periodic, being a left shift of $(\overline{i_0,\cdots,i_{N-1}},\cdots)$, depending on the $0\leq k<N$ for which $\omega\in T^kA$. It is now quite evident that the values taken by $X$ with positive probability are the $\psi_{i_k}\cdots\psi_{i_{N-1}}(c),~0\leq k<N$. 

\fin

\bigskip
\noindent
\begin{remark}
The condition of the Lemma is verified if $X(\omega)=\sum_{n\geq0}\lambda^nb(T^n\omega)$, when $b=\pm1$ and $0<\lambda<1$ is not a root of a polynomial with $0,\pm1$ as coefficients. Indeed, let $\epsilon=(\epsilon_0,\cdots,\epsilon_{n-1})$ and $\delta=(\delta_0,\cdots,\delta_{m-1})$ be minimal, with $n\leq m$. If $fix(\psi_{\epsilon})=fix(\psi_{\delta})$, then~:

$$\frac{1}{1-\lambda^n}\sum_{k=0}^{n-1}\lambda^k\epsilon_k=\frac{1}{1-\lambda^m}\sum_{k=0}^{m-1}\lambda^k\delta_k,$$

\medskip
\noindent
or $(1-\lambda^m)\sum_{k=0}^{n-1}\lambda^k\epsilon_k=(1-\lambda^n)\sum_{k=0}^{m-1}\lambda^k\delta_k$. We rewrite this as~:

\begin{eqnarray}
\sum_{k=0}^{n-1}\lambda^k(\epsilon_k-\delta_k)&=&\lambda^m\sum_{k=0}^{n-1}\lambda^k\epsilon_k-\lambda^n\sum_{k=0}^{m-1}\lambda^k\delta_k+\sum_{k=n}^{m-1}\lambda^k\delta_k\nonumber\\
&=&\({\sum_{k=n}^{m-1}\lambda^k\delta_k-\lambda^n\sum_{k=0}^{m-n-1}\lambda^k\delta_k}\)+\({\lambda^m\sum_{k=0}^{n-1}\lambda^k\epsilon_k-\lambda^n\sum_{k=m-n}^{m-1}\lambda^k\delta_k}\).\nonumber
\end{eqnarray}

\medskip
\noindent
On the right-hand side, there are only powers $\lambda$ that are $\geq n$ : between $n$ and $m-1$ in the first parenthesis and between $m$ and $n+m-1$ in the second one. As $\lambda$ is not a root of a polynomial with $0,\pm2$ coefficients, it is necessary on the left-hand side that $\epsilon_k=\delta_k$, $0\leq k<n$. Therefore $\epsilon$ is a prefix of $\delta$, which wouldn't be minimal, unless $n=m$. Thus $\epsilon=\delta$. 
\end{remark}


\providecommand{\bysame}{\leavevmode\hbox to3em{\hrulefill}\thinspace}

\bigskip
{\small{\sc{Univ Paris Est Creteil, CNRS, LAMA, F-94010 Creteil, France\\
Univ Gustave Eiffel, LAMA, F-77447 Marne-la-Vall\'ee, France}}}

\it{E-mail address~:} {\sf julien.bremont@u-pec.fr}

\end{document}